\documentclass[11pt,a4paper]{amsart}

\usepackage[utf8]{inputenc}
\usepackage[T1]{fontenc}
\usepackage{lmodern}

\usepackage{hyperref}

\usepackage{amssymb}

\theoremstyle{plain}

\newtheorem*{named_theorem}{Theorem}
\newtheorem{lemma}{Lemma}
\newtheorem*{main_proposition}{Main Proposition}
\newtheorem*{corollary}{Corollary}

\theoremstyle{definition}
\newtheorem*{definition}{Definition}

\theoremstyle{remark}
\newtheorem{remark}{Remark}

\DeclareMathOperator{\val}{V}

\newcommand{\N}{{\mathbb N}}
\newcommand{\Z}{{\mathbb Z}}

\newcommand{\Q}{{\mathbb Q}}
\newcommand{\U}{{\mathcal U}}

\newcommand{\eg}{e\textup.g\textup.{}}

\begin{document}

\title{A remark on groups without finite quotients}
\author{Jakub Gismatullin \and Alexey Muranov}
\date{\today}
\address{Instytut Matematyczny, Uniwersytetu Wroc{\l}awskiego,
pl. Grunwaldzki 2/4, 50-384 Wroc{\l}aw, Poland}
\email{gismat@math.uni.wroc.pl}
\address{Institut de Math\'ematiques de Toulouse,
Universit\'e Paul Sabatier,
31062 Toulouse Cedex 9,
France}
\email{muranov@math.univ-toulouse.fr}

\begin{abstract}
We notice that the class of nontrivial groups without proper subgroups
of finite index is not elementary, because some groups in this class,
such as $\mathbb Q*\mathbb Q$,
have ultrapowers that map homomorphically onto
$\mathbb Z/p\mathbb Z$ for every prime $p$.
Also, some ultrapowers of certain simple groups map
homomorphically onto $\mathbb Z/2\mathbb Z$.
\end{abstract}

\maketitle

\begin{definition}
By NFQ we denote the class of nontrivial groups without
proper subgroups of finite index
(equivalently, nontrivial groups which have No nontrivial Finite Quotients).
\end{definition}

For example, $(\Q,+)$ is NFQ.

Our main observation is the following proposition.

\begin{main_proposition}
If\/ $A$ and\/ $B$ are NFQ groups\textup,
then the free product\/ $G=A*B$ is NFQ as well\/\textup;
however\textup, for every non-principal ultrafilter\/ $\U$
on\/ $\omega$ and for every prime\/ $p$\textup, there exists
a homomorphism of\/ $G^{\omega}/\U$ onto\/ $\Z/p\Z$\textup,
and hence\/ $G^{\omega}/\U$ and $G^{\omega}$ are not~NFQ\textup.
\end{main_proposition}

\begin{corollary}
The class NFQ is not elementary and is not closed under
infinite products\textup.
\end{corollary}

\begin{definition}
A generating subset $S$ of a group $G$ is said to generate $G$
\emph{in\/ $n$ steps\/} if
\[
G=\underset{\text{$n$ times}}
{\underbrace{(S^{\pm1}\cup\{1\})\dotsm (S^{\pm1}\cup\{1\})}}.
\]
\end{definition}

For a group $G$ and a group word $w=w(\bar X)$, define
\[
\val_{w}(G) = \{\,w(\bar g) \,\mid\, \bar g\subset G\,\}.
\]
For example,
$\val_{X^n}(G) = \{\,g^n \,\mid\, g\in G\,\}$,
$\val_{[X,Y]}(G) = \{\,[g,h] \,\mid\, g,h\in G\,\}$.

\begin{definition}
If $w=w(\bar X)$ is a group word and $G$ a group, the
\emph{verbal width of\/ $G$ with respect to\/ $w$}
is the minimal number
of steps in which $\val_{w}(G)$ generates $\langle\val_{w}(G)\rangle$.
\end{definition}

\begin{remark}
\label{remark:1}
A group generated by its NFQ subgroups is NFQ itself.
\end{remark}

\begin{remark}
\label{remark:2}
The class NFQ is closed under taking homomorphic images, extensions,
direct sums, and free products.
\end{remark}

\begin{remark}
\label{remark:3}
An abelian group $G$ is NFQ if and only if it is divisible
(for every prime $p$, $G/pG$ is a vector space
over the finite field $\Z/p\Z$, so if it is nontrivial,
then it has an epimorphism onto $\Z/p\Z$).
\end{remark}

\begin{remark}
\label{remark:4}
An arbitrary (Cartesian) product of abelian NFQ groups is NFQ.
\end{remark}

%

\begin{remark}
\label{remark:5}
If $G$ is an NFQ group and $n\in\N$, then
$\val_{X^n}(G)\cup\val_{[X,Y]}(G)$
generates $G$
(the abelianization $G/\langle\val_{[X,Y]}(G)\rangle$ is divisible
by Remark~\ref{remark:3}).
\end{remark}

\begin{lemma}
\label{lemma:1}
Let\/ $G$ be a group and $p$ a prime number\textup.
If\/ $\val_{X^p}(G)\cup \val_{[X,Y]}(G)$
does not generate\/ $G$ in finitely many steps\textup,
then for every non-principal ultrafilter\/ $\U$ on\/ $\omega$\textup,
the ultrapower\/ $G^{\omega}/\U$ maps homomorphically
onto\/ $\Z/p\Z$\textup.
\end{lemma}

\begin{proof}
Denote $H$ the abelianization of $G^{\omega}/\U$.
Choose $f\in G^{\omega}$ such that
\[
(\forall n<\omega)\
\Bigl(f(n)\notin\val_{X^p}(G)\cdot
\underset{\text{$n$ times}}
{\underbrace{\val_{[X,Y]}(G)\dotsm \val_{[X,Y]}(G)}}\Bigr).
\]
Then $f$ represents a nontrivial element of $H/H^p$,
and hence $H/H^p$ is a nontrivial vector space over $\Z/p\Z$
and has an epimorphism onto $\Z/p\Z$.
\end{proof}

\begin{remark}
\label{remark:6}
Since every commutator is the product of $3$ squares
(\eg\ $[X,Y] = (YX)^{-2}\cdot (YX^2Y^{-1})\cdot Y^2$),
in the case $p=2$, the hypothesis of the lemma reduces to
``$\val_{X^2}(G)$ does not generate $G$ in finitely many steps.''
\end{remark}

Proof of the main proposition relies on the following remarkable
result of Rhemtulla.

\begin{named_theorem}[Rhemtulla, 1967, \cite{Rhemtulla:1968:pbefp}]
If\/ $w=w(\bar X)$ is a group word such that there exists a group\/ $H$
such that\/ $\{1\}\ne\langle\val_{w}(H)\rangle\ne H$\textup, and if\/
$A$ and\/ $B$ are two nontrivial groups of which
at least one has order greater than\/ $2$\textup,
then the verbal subgroup\/
$\langle\val_{w}(A*B)\rangle$ of\/ $A*B$
is not generated by\/ $\val_{w}(A*B)$ in finitely many steps\textup.
\end{named_theorem}

\begin{proof}[Proof of the main proposition]
Clearly $G$ is NFQ, see Remark \ref{remark:1}.

Let $w=w(X,Y,Z)=X^p[Y,Z]$.
By Rhemtulla's theorem, $G$ is not generated by $\val_{w}(G)$
in finitely many steps.
Since $\val_{w}(G)\supset\val_{X^p}(G)\cup\val_{[Y,Z]}(G)$,
$G^{\omega}/\U$ maps homomorphically onto $\Z/p\Z$ by Lemma \ref{lemma:1}.
\end{proof}

Another (more complicated)
way to prove that NFQ is not a first-order property,
without using Rhemtulla's theorem, is to consider the
simple groups constructed in \cite{Muranov:pp2009:fgisgiswvscl}:
those groups are of infinite width with respect to $X^2$,
and hence have
ultrapowers that map homomorphically onto~$\Z/2\Z$.


\end{document}